\documentclass[Svgc]{Svgc}

\usepackage{latexsym}
\usepackage{graphics}

\journalname{Graphs and Combinatorics}

\begin{document}

\title{A Blass-Sagan bijection on Eulerian equivalence classes}

\author{Beifang Chen\inst{1} \and Arthur L. B. Yang\inst{2} \and Terence Y. J.
Zhang\inst{3}}

\titlerunning{The Eulerian equivalence classes}%
\authorrunning{B. Chen et al.}%

\institute{\vskip 5pt \inst{1} Department of Mathematics, Hong Kong
University of Sciences and Technology, Clear Water Bay, Kowloon,
Hong Kong. e-mail: $^1$mabfchen@ust.hk \\[2pt] \inst{2,3} Center for Combinatorics,
LPMC-TJKLC, Nankai University, Tianjin 300071, P. R. China\\ e-mail:
$^2$yang@nankai.edu.cn, $^3$zhang@cfc.nankai.edu.cn}
\maketitle
\begin{abstract}
Following the treatment of Blass and Sagan, we present an
algorithmic bijection between the Eulerian equivalence classes of
totally cyclic orientations and the spanning trees without internal
activity edges for a given graph.
\end{abstract}


\section{Introduction}

To generalize the chromatic polynomials of graphs, Tutte
\cite{tutte1954} introduced the dichromatic polynomials in two
variables which we know as Tutte polynomials. Without much
additional effort, one can define Tutte polynomials for arbitrary
matroids. Ardila \cite{ardila2004,ardila2006} also defined the Tutte
polynomials on hyperplane arrangements. Many interesting invariants
of graphs and matroids can be computed directly from these
polynomials \cite{bjorn1992,bryox1992,crapo1967}. It is worth
mentioning that the Tutte polynomials play an important role in
statistical mechanics, where the partition functions are just simple
variants of these polynomials; the Jones polynomials and Kauffman
polynomials in knot theory are also closely related to them; see
\cite{bollo1998}. To find other new interpretations for
specializations of Tutte polynomials has interested many
mathematicians
\cite{chenbf2005,chenstanley2005,reiner1999,stanley1973,wagner1998},
etc. In this paper we concentrate on the evaluation of the Tutte
polynomial at several special points in terms of equivalence classes
of orientations on graphs.

The first remarkable result on the connection between acyclic
orientations of graphs and the Tutte polynomial is due to Stanley
\cite{stanley1973}, who gave the interpretation of the chromatic
polynomial at negative integers. Then it was generalized by Chen
\cite{chenbf2005} to interpret the integral and modular tension
polynomials of Kochol \cite{kochol2002} at nonnegative integers,
where acyclic orientations and their cut equivalence classes are
used to describe the decomposition of these polynomials. Green and
Zaslavsky \cite{grezas1983} proved that the number of acyclic
orientations with a unique source at a given vertex is the special
value of the Tutte polynomial at $(1,0)$, and a fascinating result
in \cite{chenbf2005} is that this value also counts the number of
cut equivalence classes of acyclic orientations.

Dual to Stanley's result, the number of totally cyclic orientations
also can be given by the Tutte polynomial \cite{vergna1977}.
Utilizing the theory of Ehrhart polynomials as in \cite{chenbf2005},
Chen and Stanley \cite{chenstanley2005} studied the integral and
modular flow polynomials, where they gave a similar decomposition as
the tension polynomials in terms of totally cyclic orientations and
their Eulerian equivalence classes. Dually, the number of Eulerian
equivalence classes of totally cyclic orientations is equal to the
special value of the Tutte polynomial at $(0,1)$.

Using the convolution formula due to Kook, Reiner and Stanton
\cite{krs1999}, we recover the result of Stanley in
\cite{stanley1980} which states that the value of the Tutte
polynomial at $(2,1)$ enumerates in-sequences of orientations, i.e.,
the Eulerian equivalence classes of orientations; by duality the
value at $(1,2)$ enumerates the cut equivalence classes of all
orientations. Another result from the convolution formula is the
interpretation of the value of the Tutte polynomial at $(1,1)$ in
terms of Eulerian-cut equivalence classes of orientations. Gioan
independently \cite{gioan2006} obtained the same result on
interpretations of the Tutte polynomial at $(0,1), (1,0), (1,2),
(2,1)$ and $(1,1)$, where the cycle-cocyle systems are used instead
of Eulerian-cut equivalence classes.

As Tutte originally defined, a fundamental property of the Tutte
polynomial is that it has a spanning tree expansion. Therefore,
specializations of the Tutte polynomial inherit the interpretations
in terms of spanning trees. A natural question arises: to find the
bijections between the set of some equivalence classes of
orientations and spanning trees with special property. The related
work has been done. Blass and Sagan \cite{blasag1986} constructed an
algorithmic bijection between the set of acyclic orientations and
the broken circuit complex. This algorithm was modified by Gebhard
and Sagan \cite{gebsag2000} to give a bijection between the set of
acyclic orientations with a unique sink at a given vertex and the
set of spanning trees without external activity edges. Gioan
\cite{gioan2006} gave a bijection between the set of cut equivalence
classes of acyclic orientations and the set of of acyclic
orientations with a unique sink at a given vertex. The combination
of the above two bijections leads to a bijection between cut
equivalence classes of acyclic orientations and spanning trees
without external activity edges. Gioan and Vergnas \cite{giover2005}
also established the activity preserving bijections between spanning
trees and orientations.

The main task of this paper is to give a Blass-Sagan bijection
between Eulerian equivalence classes of totally cyclic orientations
and spanning trees without internal activity edges. As each cut
equivalence class of acyclic orientations has an acyclic orientation
with a unique sink at a given vertex, our bijection would be helpful
to find the corresponding representative element for each Eulerian
equivalence class of totally cyclic orientations.

\section{Definitions and notations}

Let $G=(V,E)$ be a graph with vertex set $V$ and edge set $E$, in
which  multiple edges and loops are allowed. Given $e\in E$, let
$G-e=(V,E\backslash \{e\})$. Thus $G-e$ is obtained from $G$ by
deleting the edge $e$. Let $G/e$ be the multigraph obtained from $G$
by contracting the edge $e$. Throughout this paper the graphs are
assumed to be always connected.

Now let us define the Tutte polynomial $T_G(x,y)$ for a graph $G$
recursively. First, let $T_{E_n}(x,y)=1$, where $E_n$ is the empty
$n$-graph for $n\geq 1$. In general, we have
$$T_G(x,y)=\left\{
\begin{array}{ll}
xT_{G/e}(x,y) & \mbox{if $e$ is a bridge,}\\
yT_{G-e}(x,y) & \mbox{if $e$ is a loop,}\\
T_{G-e}(x,y)+T_{G/e}(x,y) & \mbox{if $e$ is neither a bridge nor a
loop.}
\end{array}
\right.
$$

As we remarked at the beginning, the original definition of
$T_G(x,y)$ is in terms of spanning trees of $G$. We adopt the
notions of \cite{bollo1998} in the following. For a connected graph
$G=(V,E)$, a tree $F=(V',E')$ is a spanning tree of $G$ if $V'=V$
and $E'\subset E$. If $G$ is not connected, the spanning trees of
all components form a spanning forest of $G$. Now let us impose an
order on the edge set $E(G)=\{e_1, e_2, \ldots, e_m\}$, with $e_i$
preceding $e_j$ if $i<j$. Fix a spanning forest $F$ of $G$. For each
edge $e_i$ in $F$, we call $U_F(e_i)=\{e_j\in E(G) :
(F-e_i)+e_j\mbox{ is a spanning forest}\}$ the cut defined by $e_i$.
If $e_i$ is the smallest edge of the cut it defines, we call $e_i$
an internally active edge of $F$. Similarly, for each edge $e_j$ not
in $F$, we call $Z_F(e_j)=\{e_i\in E(G): e_i \mbox{ is an edge on
the unique cycle of } F+e_j\}$ the cycle defined by $e_j$. If $e_j$
is the smallest edge of the cycle it defines, we call $e_j$ an
externally active edge. We say that a spanning forest has internal
activity $i$ and external activity $j$ if there are precisely $i$
internally active edges and precisely $j$ externally active edges,
denoted by an $(i,j)$-forest. Tutte originally defined
\begin{equation}\label{tutte-def}
T_G(x,y)=\sum_{i,j}t_{ij}x^iy^j,
\end{equation}
 where $t_{ij}$ is the number of
$(i,j)$-forests.

Recall that a cut of $G$ is a partition $[S, T]$ of the vertex set
$V$ such that the removal of $[S, T]$, the set of all edges between
$S$ and $T$, disconnects the graph $G$. For a digraph $(G,
\varepsilon)$, where $\varepsilon$ is an orientation of $G$, we
denote by $(S, T)_\varepsilon$ the set of all edges going from $S$
to $T$, and by $(T, S)_\varepsilon$ the set of all edges going from
$T$ to $S$. A bond is a minimal cut. A bond $[S,T]$ is called
directed relative to $\varepsilon$ if $(S, T)_\varepsilon=\emptyset$
or $(T, S)_\varepsilon=\emptyset$. A cut is called directed if it
can be decomposed into a disjoint union of directed bonds. Let
$\mathcal{O}(G)$ denote the set of all orientations of G,
$\mathcal{AO}(G)$ the set of all orientations without directed
cycles, and $\mathcal{BO}(G)$ the set of all orientations without
directed cuts.

Given an orientation $\varepsilon$ of $G$, a directed edge $e=(u,v)$
is called cut flippable if there are no directed paths either from
$u$ to $v$ or from $v$ to $u$ in $G-e$. An directed edge $e$
relative to $\varepsilon$ is called cycle flippable if there are
directed paths both from $u$ to $v$ and  from $v$ to $u$ in $G-e$.

We call two orientations $\varepsilon_1$ and $\varepsilon_2$
cut-equivalent, denoted by $\varepsilon_1\sim_{c} \varepsilon_2$, if
the spanning subgraph induced by the edge set $\{e\in E(G)\ |\
\varepsilon_1(e)\neq \varepsilon_2(e)\}$ is a directed cut with
respect to $\varepsilon_1$ or $\varepsilon_2$. It is easy to see
that $\sim_{c}$ is an equivalence relation on $\mathcal{O}(G)$, and
it also induces an equivalence relation on $\mathcal{AO}(G)$.

Similarly, we define the Eulerian equivalence relations as follows.
We call two orientations $\varepsilon_1$ and $\varepsilon_2$
Eulerian equivalent, denoted by $\varepsilon_1\sim_{e}
\varepsilon_2$, if the spanning subgraph induced by the edge set
$\{e\in E(G)\ |\ \varepsilon_1(e)\neq \varepsilon_2(e)\}$ is a
directed Eulerian graph with respect to $\varepsilon_1$ or
$\varepsilon_2$, i.e., the in-degree is equal to the out-degree at
each vertex. It is easy to see that $\sim_{e}$ is an equivalence
relation on $\mathcal{O}(G)$, and it induces an equivalence relation
on $\mathcal{BO}(G)$.

We also need the concept of Eulerian-Cut equivalence over
orientations. Two orientations $\varepsilon_1$ and $\varepsilon_2$
are called to be Eulerian-cut equivalent, denoted by
$\varepsilon_1\sim_{ec} \varepsilon_2$, if the spanning subgraph
induced by the edge set $\{e\in E(G)\ |\ \varepsilon_1(e)\neq
\varepsilon_2(e)\}$ is a disjoint union of a directed Eulerian graph
and a direct cut with respect to $\varepsilon_1$ or $\varepsilon_2$.
The relation $\sim_{ec}$ is also an equivalence relation on
$\mathcal{O}(G)$.

\begin{figure}
\begin{center}
  \includegraphics{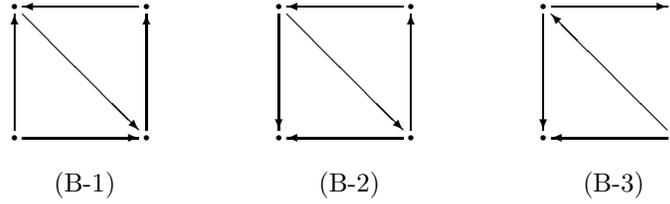}
\caption{Equivalence relations among three orientations.}
\label{fig-1}
\end{center}
\end{figure}

By definitions, the two orientations (B-1) and (B-2)  in Fig.
\ref{fig-1} are cut equivalent, (B-2) and (B-3) are Eulerian
equivalent, while (B-1) and (B-3) are Eulerian-cut equivalent.

\section{Eulerian equivalence classes}

Using the theory of Ehrhart polynomials, Chen and Stanley obtained
the following nice result, which is independently discovered by
Gioan \cite{gioan2006}.

\begin{theorem}\cite[Theorem 1.2]{chenstanley2005} \label{ec-thm}
For any graph $G$, let $\alpha(G)$ denote the number of Eulerian
equivalence classes of $\mathcal{BO}(G)$. Then
\begin{equation}
\alpha(G)= T_G(0,1).
\end{equation}
\end{theorem}

In the following we will present two proofs of the above theorem.
The first proof is purely inductive according to the inductive
definition of Tutte polynomials, and the second one is an
algorithmic bijection similar to the modified Blass-Sagan algorithm
\cite{gebsag2000}.

\subsection{The inductive proof}

 For any fixed edge $e=(u,v)$, it is clear that there always exists an orientation
 $\varepsilon$ in each Eulerian equivalence class of
 $\mathcal{BO}(G)$ such that the edge $e$ is directed from $u$ to
 $v$ with respect to $\varepsilon$. Notice that the edge $e$ has the same cycle flippable property in each Eulerian
 equivalence class, i.e., for any two equivalent totally cyclic orientations $\varepsilon$ and  $\varepsilon'$ with $\varepsilon(e)=\varepsilon'(e)$,
 then $e$ is cycle flippable relative to $\varepsilon$ if and only
 if it is cycle flippable relative to $\varepsilon'$. Therefore,
 in each equivalence class we can choose an orientation with $e$ directed from $u$ to
 $v$ as a representative element.

\noindent \textit{Proof of Theorem \ref{ec-thm}.}  We shall deduce
the assertion from the following four properties of the function
$\alpha(G)$.

\begin{itemize}
\item[(i)] If $G=E_n$, then $\alpha(G)=1$.

\item[(ii)] If $e$ is a loop, then $\alpha(G)=\alpha(G-e)$.

\item[(iii)] If $e$ is a bridge, then $G$ has no totally cyclic orientations so $\alpha(G)=0$.

\item[(iv)] Finally, suppose that $e$ is neither a bridge nor a
loop. Consider an equivalence class of $\mathcal{BO}(G)$, and the
orientation $\varepsilon$ is its representative element. If $e$ is
cycle flippable relative to $\varepsilon$, then all orientations
equivalent to $\varepsilon$ give an equivalence class of
$\mathcal{BO}(G-e)$; otherwise, they give an equivalence class of
$\mathcal{BO}(G/e)$. Also, all appropriate equivalence classes of
$\mathcal{BO}(G-e)$ and $\mathcal{BO}(G/e)$ arise in this way.
Therefore, in this case we have
$$\alpha(G)=\alpha(G-e)+\alpha(G/e).$$
\end{itemize}
Since $\alpha(G)$ and $T_G(0,1)$ satisfy the same boundary
conditions and recurrence relations, the desired result immediately
follows. \qed

\subsection{The bijective proof}

From Equation (\ref{tutte-def}) we see that the value $T_G(0,1)$
counts the number of spanning trees without internal activity edges.
To prove Theorem \ref{ec-thm}, it suffices to establish a bijection
between these spanning trees of $G$ and Eulerian equivalence classes
of $\mathcal{BO}(G)$.

Fix an orientation $\varepsilon$ of $G$ (not necessarily totally
cyclic or acyclic), which we will refer to as the {\em  normal
orientation}. Fix the total order imposed on the edges which defines
the internal and external activity. We say that an orientation
$\varepsilon'$ is {\em  reduced} if for each edge $e\in E(G)$ either
$\varepsilon(e)=\varepsilon'(e)$ or there exists no directed cycle
containing $e$ with other edges smaller than $e$.

For any oriented arc $e=\vec{uv}$, we denote the oppositely oriented
arc by $e'=\vec{vu}$. To {\em  unorient} an arc $e$ for an
orientation $\varepsilon$ of $G$, it means that we will just add the
oppositely oriented arc $e'$. Given a graph with unoriented edges,
let $G'$ be the {\em  contraction} of $G$, which is the graph where
all unoriented edges have been contracted. The orientation of $G'$
is inherited from the original graph $G$. We say that $G$ is reduced
if its contraction $G'$ is reduced with respect to the inherited
normal orientation. For any two orientations $\varepsilon_1$ and
$\varepsilon_2$ of $G$ with unoriented edges, we say that they are
Eulerian equivalent if the two inherited orientations of the
contraction $G'$ are Eulerian equivalent.

\begin{lemma} For the normal orientation $\varepsilon$ and the total order on edges fixed as above,
there exists one and only one reduced orientation in each Eulerian
equivalence class of $\mathcal{BO}(G)$.
\end{lemma}

\proof Given an Eulerian equivalence class, we first show that there
exists at least one reduced orientation. Start with one arbitrary
totally cyclic orientation, say $\varepsilon_0$. If $\varepsilon_0$
is reduced, then we are done. Otherwise, find the largest edge, say
$e_0$, which doesn't satisfy the reduced property. It means that
$\varepsilon(e_0)\neq \varepsilon_0(e_0)$ and there exists one
directed cycle which contains $e_0$ and all other edges on the cycle
are smaller than $e_0$. By reversing the orientation of this cycle,
we obtain another Eulerian equivalent orientation $\varepsilon_1$
with all edges larger than or equal to $e_0$ satisfying the reduced
property. Iterating the above process, we will get one orientation
equivalent to $\varepsilon_0$, with all its edges satisfying the
reduced property.

Now we show that the reduced orientation is unique in the Eulerian
equivalence class. Suppose there are two reduced equivalent
orientations $\varepsilon'$ and $\varepsilon''$. Consider the
spanning subgraph induced by the edge set $\{e\in E(G)\ |\
\varepsilon'(e)\neq \varepsilon''(e)\}$. If not empty, then it must
contain a directed cycle with respect to $\varepsilon'$ or
$\varepsilon''$. Therefore, the largest edge on this cycle satisfies
the reduced property only for one of two orientations $\varepsilon'$
and $\varepsilon''$. This is a contradiction. \qed

As shown above, from an arbitrary orientation $\varepsilon'$ we can
obtain the reduced orientation in each Eulerian equivalence class
with the iterated process. For convenience we call it the {\em
normalization} of $\varepsilon'$.

In the following we will construct an algorithm which maps each
reduced totally cyclic orientation to a spanning tree without
internal activity edges. Due to the above lemma, we obtain the
desired bijection. With the total order imposed on the edge set,
each oriented edge is sequentially examined and is either  deleted
or unoriented using the following algorithm:

\begin{itemize}

\item[(S1)] Input a graph $(G, \varepsilon)$, where $\varepsilon$
is an orientation of $G$ with some unoriented edges.

\item[(S2)] Let $(G',\varepsilon')$ be the contraction of $(G,
\varepsilon)$ with all unoriented edges having been contracted. If
$\varepsilon'$ is not reduced, then we take the reduced
representation $\varepsilon''$ in its Eulerian equivalence class.

\item[(S3)] Consider the largest edge $e$ of $G'$. If $e$ is  a
loop or cycle flippable with respect to $\varepsilon''$, then we
delete $e$ from $G'$. Otherwise, we unorient $e$ in $G'$. Reset $G$
to be the graph recovered from $G'$ by adding back all unoriented
edges. Reset $\varepsilon$ to be the orientation of $G$ obtained
from $\varepsilon''$, i.e., for all oriented edge $e'$ we have
$\varepsilon''(e')=\varepsilon(e')$. If $G$ contains at least one
oriented edge with respect to $\varepsilon$, then go to Step (S2).
Otherwise, go to Step (S4).

\item[(S4)] Output the graph $G$.

\end{itemize}

\begin{figure}
\begin{center}
  \includegraphics{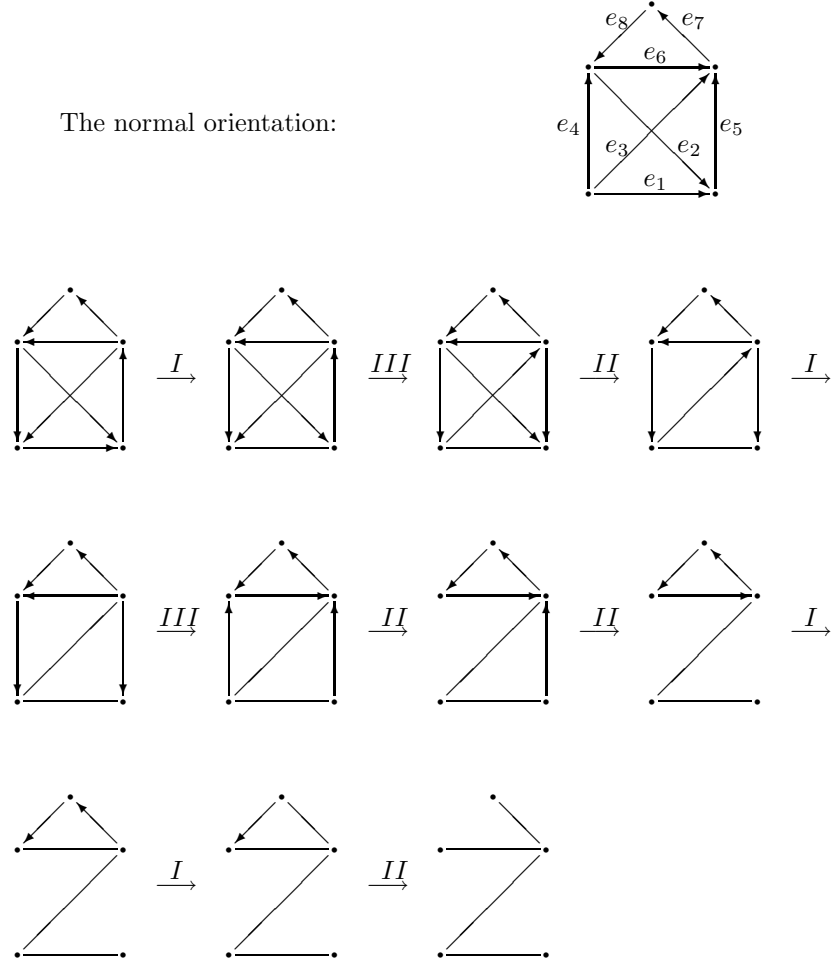}
\caption{An example of the algorithm} \label{fig-alg}
\end{center}
\end{figure}

For an example of how the above algorithm works, see Figure
\ref{fig-alg}, where $I$ denotes the unorientation, $II$ denotes the
deletion, and $III$ denotes the normalization.

To show that this algorithm actually does produce a bijection, we
shall first introduce a sequence of sets, $\mathcal{O}_0,
\mathcal{O}_1, \cdots, \mathcal{O}_q$, such that $\mathcal{O}_0$ is
the set of all reduced totally cyclic orientations of $G$, and
$\mathcal{O}_q$ (where $q=|E(G)|$) is the set of all spanning trees
of $G$ without internal activity edges. We will show that the $k$-th
step of the algorithm gives a bijection, $f_k:
\mathcal{O}_{k-1}\rightarrow \mathcal{O}_k $, where $\mathcal{O}_k$
is defined as the set of all orientations $\varepsilon$ of spanning
subgraphs of $G$ satisfying the following conditions:
\begin{itemize}

\item[(a)] Each of the first $k$ largest edges of $G$ is either present in
$\varepsilon$ (as an unoriented edge) or absent from $\varepsilon$,
but each of the remaining $q-k$ edges is present in $\varepsilon$ in
exactly one orientation, and there does not exist a cycle only
consisting of unoriented edges.

\item[(b)] $\varepsilon$ is totally cyclic.

\item[(c)] $\varepsilon$ is reduced.

\item[(d)] For each unoriented edge $e$ in the subgraph, if $e$ is a
bridge which separates the subgraph into two components $C_1$ and
$C_2$, there exists at least one edge strictly smaller than $e$ in
the edge cut $E_G[C_1,C_2]$.

\end{itemize}

From the above conditions, we see that $\mathcal{O}_0$ is indeed the
set of all reduced representations of the totally cyclic
orientations of $G$, and $\mathcal{O}_q$ is indeed the set of all
spanning trees without internal activity edges.

\begin{lemma}
$f_k$ maps $\mathcal{O}_{k-1}$ into $\mathcal{O}_k$.
\end{lemma}

\proof It suffices to verify that properties (a)--(d) listed
previously are still satisfied after the algorithm is applied at the
$k$-th stage.

\begin{itemize}
\item[(a)] If the $k$-th largest edge $e$ is cycle flippable then the algorithm will delete it; otherwise, the algorithm will unorient
it. Therefore, it will not create a new cycle consisting of only
unoriented edges.

\item[(b)] Clearly, to unorient an edge and delete the cycle flippable edge will not
destroy the totally cyclic property.

\item[(c)] This is ensured by Step (S2) of the algorithm.

\item[(d)] Suppose that there exists some unoriented  edge $e$  as a bridge
in the subgraph such that  $e$ is the smallest edge in the edge cut
$E_G[C_1,C_2]$. Therefore, in the process of the algorithm all edges
of $E_G[C_1,C_2]$ except $e$ will be deleted, i.e., all these
oriented edges are cycle flippable. Now consider the second smallest
edge $e_0$ of $E_G[C_1,C_2]$. Clearly, $e_0$ must not be cycle
flippable, and the algorithm will unorient it. This is a
contradiction.
\end{itemize}
\qed

To prove that $f_k$ is bijective, we first give the following two
lemmas:

\begin{lemma}\label{lem-del}
Given an orientation $\varepsilon\in\mathcal{O}_{k-1}$, let $e$ be
the largest oriented edge of the underlying graph $G$. Let $G'=G-e$,
and $\varepsilon'$ be the orientation of $G'$ inherited from
$\varepsilon$. If $\varepsilon$ is reduced and $e$ is cycle
flippable in $\varepsilon$, then $\varepsilon'$ is reduced.
Moreover, $f_k(\varepsilon(G))=\varepsilon'(G')$.
\end{lemma}

\proof Suppose that $\varepsilon'$ is not reduced. There must exist
one edge $e'$ which is smaller than $e$ and doesn't satisfy the
reduced property in $G'$. Clearly, $e'$ also doesn't satisfy the
reduced property for the orientation $\varepsilon$ in $G$, which is
contrary to the fact that $\varepsilon$ is reduced. \qed

\begin{lemma} \label{lem-uno} Given any two distinct reduced totally cyclic orientations $\varepsilon_1$ and
$\varepsilon_2$ of $G$, suppose that the largest oriented edge $e$
is neither cycle flippable with respect to $\varepsilon_1$ nor
$\varepsilon_2$. Let $\varepsilon_1'$ (resp. $\varepsilon_2'$) be
the orientation of $G$ obtained from $\varepsilon_1$ (resp.
$\varepsilon_2$) by unorienting the edge $e$. Then $\varepsilon_1'$
and $\varepsilon_2'$ are not Eulerian equivalent as orientations of
the contraction graph $G/e$.
\end{lemma}

\proof Since $\varepsilon_1, \varepsilon_2$ are reduced and $e$ is
the largest edge in $G$, we must have $\varepsilon_1(e)=
\varepsilon_2(e)$. Suppose that $\varepsilon_1'$ and
$\varepsilon_2'$ are Eulerian equivalent, then the edge set $\{e'\in
E(G/e)\ |\ \varepsilon_1'(e')\neq \varepsilon_2'(e')\}$ can be taken
as a disjoint union $\cup_i C_i$, where each $C_i$ is a directed
cycle in $G/e$ with respect to $\varepsilon_1'$ or $\varepsilon_2'$.
The set $\{e'\in E(G/e)\ |\ \varepsilon_1'(e')\neq
\varepsilon_2'(e')\}$ can not be empty, otherwise we will have
$\varepsilon_1 \sim_e \varepsilon_2$, contradicting the fact that
they are distinct reduced orientations. If for each $i$ the edges in
$G$ corresponding to the edges of $C_i$ also form a cycle, then we
also have $\varepsilon_1 \sim_e \varepsilon_2$. Otherwise, suppose
for some $i$ the edges in $G$ corresponding to the edges of $C_i$ do
not form a cycle, but together with the edge $e$ they will form a
cycle. If $C_i$ and $e$ form a directed cycle with respect to
$\varepsilon_1$ (resp. $\varepsilon_2$), then $e$ will be cycle
flippable with respect to $\varepsilon_2$ (resp. $\varepsilon_1$),
which is again a contradiction. \qed

\begin{theorem}  $f_k$ is bijective.
\end{theorem}
\proof First we prove that $f_k$ is one to one. Suppose
$\varepsilon_1$ and $\varepsilon_2$ are two distinct elements of
$\mathcal{O}_{k-1}$ which are both mapped to $\varepsilon\in
\mathcal{O}_{k}$ by the algorithm. Since the algorithm only affects
the $k$-th large edge, we note that in both $\varepsilon_1$ and
$\varepsilon_2$, the cases are same for the first $k-1$ large edges
of $G$. We note that $\varepsilon$ was not obtained from
$\varepsilon_1$ and $\varepsilon_2$ by deletion. Otherwise,
$\varepsilon_1$ and $\varepsilon_2$ will be the same due to Lemma
\ref{lem-del}. Thus we only need consider the case that
$\varepsilon$ was obtained from $\varepsilon_1$ and $\varepsilon_2$
by unorienting the $k$-th edge and applying the normalization. By
Lemma \ref{lem-uno}, this is also impossible.

Then we prove that that $f_k$ maps $\mathcal{O}_{k-1}$ onto
$\mathcal{O}_k$. For any $\varepsilon\in\mathcal{O}_k$ such that the
$k$-th edge $e$ of $G$ is absent in the underlying spanning
subgraph, we just add the edge $e$ in the subgraph and normally
orient it. Denote the orientation of this new diagraph by
$\varepsilon'$. Since $\varepsilon$ is totally cyclic and the
underlying graph is connected, $\varepsilon'$ is still totally
cyclic. Notice that $e$ is the largest oriented edge with respect to
$\varepsilon'$. Therefore, $\varepsilon'$ is also reduced and the
directed edge $e$ is cycle flippable. It means that
$\varepsilon'\in\mathcal{O}_{k-1}$, and the $k$-th stage of the
algorithm will map $\varepsilon'$ to $\varepsilon$.

For any $\varepsilon\in\mathcal{O}_k$ such that the $k$-th edge $e$
of $G$ is unoriented in the underlying spanning subgraph, we
construct one orientation $\varepsilon'\in\mathcal{O}_{k-1}$ as
follows.

\begin{itemize}
\item[(1)] Choose an orientation of $e$
such that the new orientation is totally cyclic. Note that such an
orientation always exists.

\item[(2)] Normalize the new orientation. If the directed edge $e$ is not cycle flippable, then return
the orientation. Otherwise, go to (3).

\item[(3)] Reorient the edge $e$
oppositely, then reorient the directed cycle containing $e$
oppositely, and go to (2).
\end{itemize}

Let $\varepsilon'$ be the returned orientation. Clearly, $e$ is not
cycle flippable  with respect to $\varepsilon'$, and
$\varepsilon'\in\mathcal{O}_{k-1}$. The $k$-th stage of the
algorithm will map $\varepsilon'$ to $\varepsilon$. \qed

\begin{remark}
In fact, the acyclic orientations with only one given source (or
sink) can be considered as the representative elements of cut
equivalence classes of acyclic orientations. But for Eulerian
equivalence classes of totally cyclic orientations, the dual
representative elements are not known. E. Gioan mentioned to use the
degree sequences to characterize the Eulerian equivalence classes.
In this paper our reduced orientations are also representative
elements of Eulerian equivalence classes, but they depend on the
total order on the edge set and the fixed normal orientation.
\end{remark}

\begin{acknowledgement}
This work was supported by the 973 Project on Mathematical
Mechanization, the PCSIRT Project of the Ministry of Education, the
Ministry of Science and Technology, and the National Science
Foundation of China. The second author would like to thank Professor
Beifang Chen for his hospitality during the visit to HKUST.
\end{acknowledgement}

\end{document}